\documentclass[12pt]{amsart}
\usepackage{amscd,times,fullpage}
\usepackage[utf8x]{inputenc}

\usepackage{amsthm,amsmath}
\usepackage{amssymb}
\usepackage{graphicx}
\usepackage{amsthm}
\usepackage{enumerate}
\usepackage{dsfont}
\usepackage{mathtools}
\usepackage{tikz-cd}
\usetikzlibrary{cd}
\usetikzlibrary{positioning}
\usetikzlibrary{matrix}
\usepackage[hidelinks]{hyperref}

\newcommand{\R}{\mathbb{R}}

\newcommand{\Q}{\mathbb{Q}}
\newcommand{\CP}{\mathbb{C}\mathrm{P}}

\newcommand{\CH}{\mathbb{C}\mathrm{H}}

\newcommand{\Ric}{\mathrm{Ric}}

\newcommand{\C}{\mathbb{C}}            
\newcommand{\de}{\partial}

\newcommand{\K}{K\"{a}hler }

\newcommand{\F}{\mathcal F}

\newcommand{\ov}[1]{\overline{#1}}

\newcommand{\di}{{\operatorname{d}}}

\newcommand{\Id}{\operatorname{Id}}

\newcommand{\lra}{\longrightarrow}

\newcommand{\D}{\mathcal{D}}

\DeclarePairedDelimiter{\norm}{\lVert}{\rVert}
\DeclarePairedDelimiter{\abs}{\lvert}{\rvert}
\newcommand{\CO}{{{\mathcal C} \times \Omega}}
\newcommand{\hs}{\hspace{0.1em}}
\newcommand{\FR}{{\mathfrak{R}}}
\newcommand{\FF}{{\mathfrak{F}}}

\newcommand{\Lie}{\mathcal{L}}

\newtheorem{theor}{Theorem}
\newtheorem{prop}{Proposition}

\newtheorem{lem}[prop]{Lemma}
\newtheorem{cor}{Corollary}

\newtheorem{ex}{Example}
\newtheorem{remark}[prop]{Remark}

\newcommand{\mc}{{\mathcal C}}
\newcommand{\me}{{\mathcal E}}

\newcommand{\g}{{\operatorname {cig}}}

\makeatletter
\newtheorem*{rep@theorem}{\rep@title}
\newcommand{\newreptheor}[2]{%
	\newenvironment{rep#1}[1]{%
		\def\rep@title{#2 \ref{##1}}%
		\begin{rep@theorem}}%
		{\end{rep@theorem}}}
\makeatother

\newreptheor{theor}{Theorem}

\usepackage{nicefrac}

\begin{document}

\title[Immersions of Sasaki-Ricci solitons into homogeneous Sasakian manifolds]{Immersions of Sasaki-Ricci solitons into homogeneous Sasakian manifolds}
	
\author{R.~Mossa}
\address{Dipartimento di Matematica e Informatica, Universit\'a degli studi di Cagliari, Via Ospedale 72, 09124 Cagliari, Italy}
\email{roberto.mossa@unica.it}
\author{G.~Placini}
\address{Dipartimento di Matematica e Informatica, Universit\'a degli studi di Cagliari, Via Ospedale 72, 09124 Cagliari, Italy}
\email{giovanni.placini@unica.it}
{

\date{\today ; {\copyright  R.~Mossa and G.~Placini 2023}}

\thanks{The authors are supported by INdAM and  GNSAGA - Gruppo Nazionale per le Strutture Algebriche, Geometriche e le loro Applicazioni and by GOACT - Funded by Fondazione di Sardegna. The second author is  funded by the National Recovery and Resilience Plan (NRRP), Mission 4 Component 2 Investment 1.5 - Call for tender No.3277 published on December 30, 2021 by the Italian Ministry of University and Research (MUR) funded by the European Union – NextGenerationEU. Project Code ECS0000038 – Project Title eINS Ecosystem of Innovation for Next Generation Sardinia – CUP F53C22000430001- Grant Assignment Decree No. 1056 adopted on June 23, 2022 by (MUR)}

\begin{abstract}
We discuss local Sasakian immersion of Sasaki-Ricci solitons (SRS) into fiber products of homogeneous Sasakian manifolds. In particular, we prove that SRS locally induced by a large class of fiber products of homogeneous Sasakian manifolds are, in fact, $\eta$-Einstein. The results are stronger for immersions into Sasakian space forms.
Moreover, we show an example of a K\"ahler-Ricci soliton on $\mathbb C^n$ which admits no local holomorphic isometry into products of homogeneous bounded domains with flat K\"ahler manifolds and generalized flag manifolds.
\end{abstract}

\maketitle

\section{Introduction and statements of the main results} \label{sectionint}

Sasaki-Ricci solitons (SRS for short) were introduced in \cite{FutakiOnoWang09DefSRS} as special solutions of the Sasaki-Ricci flow of \cite{Smoczyk10SRFlowDefinition}. Since then the Sasaki-Ricci flow and SRS have received growing interest, see for instance \cite{Collins16SRFlowConvergenceCRELLE,Collins15SRFlowConvergence,Petrecca16SRSDeformations,Tadano15GapSRS} and references therein.
Given the analogy between Sasakian and \K geometry, it is no surprise that SRS are closely related to K\"ahler-Ricci solitons (KRS for short). 
Roughly speaking, Sasaki-Ricci solitons can be thought of as KRS transverse to the Reeb foliation.
Motivated by recent results in the \K setting (\cite{LoiMossa23Rigidity}) we investigate SRS induced by  immersions into Sasakian space forms and, more generally, homogeneous Sasakian manifolds.
Our first result deals with Sasaki-Ricci solitons $M$ locally induced by a non-elliptic Sasakian space form $S(N,c)$. In this case we show that $M$ is necessarily a quotient of a Sasakian space form of the appropriate dimension and therefore the SRS is trivial, i.e., $M$ is an $\eta$-Einstein manifold.

\begin{theor}\label{ThmSasakiSolitonNoncompact}
	Let  $M$  be a $(2n+1)$-dimensional complete Sasakian manifold endowed with a Sasaki-Ricci soliton. 
	Suppose there exist a neighbourhood $U_p$ of a point $p\in M$ and an immersion $\psi\colon U_p\lra S(N,c)$ into a Sasakian space form $S(N,c)$ with $c\leq -3$ .
	
	Then $M$ is Sasaki equivalent to $S(n,c)/\Gamma$ for a discrete subgroup $\Gamma$ of the group of Sasakian transformations of $S(n,c)$. In particular, $M$ is $\eta$-Einstein.
	
	If additionally $U_p=M$, then $\Gamma=0$ and, up to a Sasakian transformation of $S(n,c)$, $\psi$ is of the form $$\psi(z,t)=(z,0,t+a)$$ for a constant $a$.
\end{theor}
The behaviour is less rigid in the elliptic case. 
For instance, all compact homogeneous Sasakian manifolds can be realized as $\eta$-Einstein Sasakian submanifolds of the standard Sasakian sphere, see \cite[Theorem~1.4]{LoiPlaciniZedda23SasakiHomogeneous}.
Therefore, one cannot hope to get an equally strong conclusion in the compact case.
Nevertheless, we show that SRS locally induced by the standard Sasakian sphere are necessarily $\eta$-Einstein with rational constants.

\begin{theor}\label{ThmSasakiSolitonCompact}
	Let  $M$  be a $(2n+1)$-dimensional complete Sasakian manifold endowed with a Sasaki-Ricci soliton. 
	Suppose there exists a neighbourhood $U_p$ of a point $p\in M$ and an immersion $\psi\colon U_p\lra S^{2N+1}$ into the standard Sasakian sphere.
	Then $(M,g)$ is a $\eta$-Einstein Sasakian manifold whose $\eta$-Einstein constants $(\lambda,\nu)$ are given by $\lambda=4\mu -2$ for some $\mu\in\Q$.
\end{theor} 
Under the additional hypothesis of regularity the same conclusions of Theorems~\ref{ThmSasakiSolitonNoncompact} and \ref{ThmSasakiSolitonCompact} were reached in \cite{Placini21SasakiSolitons} by considering the holomorphic isometries of the base of the Boothby-Wang fibration of $M$ (which only exists because $M$ is regular).
We are able to get rid of  the regularity assumption by combining the proof in \cite{Placini21SasakiSolitons} with the results of \cite{LoiPlaciniZedda23SasakiHomogeneous} on the extension of immersions into Sasakian space forms.

As mentioned above, up to transverse homotheties, all homogeneous Sasakian manifolds can be realized as Sasakian submanifolds of the standard Sasakian sphere. Thus we get the following
immediate consequence of Theorem~\ref{ThmSasakiSolitonCompact}.

\begin{cor}\label{CorHomogeneous}
	Let $M$ be a Sasakian manifold endowed with a SRS which admits a local immersion $\varphi:U\lra N$ at a point $p\in M$ in a compact Sasakian homogeneous manifold $N$. 
	Then $M$ is $\eta$-Einstein. In particular, on a compact manifold $N$ there exist no non-trivial homogeneous Sasaki-Ricci solitons.
\end{cor}
In this case we cannot deduce the rationality of the Einstein constants of the metric becuse this depends on the homothety that is performed in order to embed the homogeneous Sasakian manifold $N$ into a standard sphere. Namely, the Einstein constants are rational if and only if the factor of the transverse homothety is rational.

It is natural to ask whether similar results hold if one replaces Sasakian space forms with homogeneous Sasakian manifolds. The analogous question was investigated and partially answered in the recent works of the first named author and Loi \cite{LoiMossa21KRSIntoComplexSpaceForms,LoiMossa22HBD,LoiMossa22BlowUp,LoiMossa23Rigidity} in the \K case. 
In full generality the answer turns out to be negative both in the \K and in the Sasakian setting, see Examples~\ref{ExSoliton1} and \ref{ExSoliton2} in Section~\ref{sectionbackground}.
Nevertheless, we can answer the question in the affirmative if we specialize to specific model manifolds.
Namely, one can easily construct Sasakian manifolds as Sasakian fiber products of Sasakian space forms (or, more generally, regular Sasakian manifolds) by identifying the leaves of the Reeb foliation, see Section~\ref{sectionbackground} below.
In analogy with the \K case, we consider SRS that can be locally immersed into the Sasakian fiber product of a compact homogeneous Sasakian manifold $N$ and a Sasakian manifold that fibers over a homogeneous bounded domain $(\Omega,g_\Omega)$.
\begin{theor}\label{ThmProducts}
	Let $M$ be a manifold endowed with a SRS and let $M'=\Omega\times N$ be the product of a homogeneous bounded domain in $\C^n$ and a compact homogeneous Sasakian manifold $N$ endowed with the standard Sasakian product structure. 
	If there exists a local Sasakian immersion $\psi:U_p\lra M'$ for a point $p\in M$, then $M$ is $\eta$-Einstein.
\end{theor}
\begin{remark}\rm
	Observe that each factor of a Sasakian fiber product can be naturally realized as a Sasakian submanifold. Therefore, Theorem~\ref{ThmProducts} implies that non-trivial SRS $M$ cannot be locally induced by immersions $M\lra \Omega\times \R$ or $M\lra N$. 
	Notice that, while the latter is prevented by Corollary~\ref{CorHomogeneous}, the former does not follow from other known results.
\end{remark}

This result can be considered the Sasakian analogue of \cite[Theorem~1.1]{LoiMossa23Rigidity}
The choice of the model manifold $\Omega\times N$ is not arbitrary.
In fact, as shown in Examples~\ref{ExSoliton1} and \ref{ExSoliton2}, there exist SRS on the Sasakian fiber product of $\R^3$ and another Sasakian space forms.
Therefore, given a SRS on a manifold $M$, it is not a priori clear whether this can be locally induced by a Sasakian fiber product such that one factor is the Sasakian space form $\R^{2N+1}$.
Motivated by this observation, we study immersions  into the Sasakian fiber product $\C^m\times N\times\Omega$ where $\C^m$ is endowed with the flat metric, $N$ is a compact homogeneous Sasakian manifold and $\Omega$ is a homogeneous bounded domain with a homogeneous metric $g_\Omega$.
The SRS we consider is the  cigar metric on $\R^{2n+1}$ induced by the cigar soliton $g_\g$ on $\C$.
Namely, consider the cigar metric $g_\g$ on $\C$ presented by Hamilton in \cite{Hamilton88RicciFlow}. This gives rise to a simple SRS on $R^3$. 
Moreover, by taking Sasakian fiber products of copies of $\R^3$, it provides a SRS on $\R^{2n+1}$ for all $n$, see Example~\ref{ExSasakianCigar} below.
By abuse of terminology, we call this metric on $\R^{2n+1}$ the cigar metric.
A local Sasakian immersion of $\R^3$ with the cigar metric into $\C^m\times N\times\Omega$ covers a holomorphic isometry of $\left(\C, g_{\g}\right)$ into $\C^m\times \mathcal{C}\times\Omega$ where $\mathcal{C}$ is a generalized flag manifold of integral type.
Since a large enough multiple of the \K metrics of the three factors admits an immersion into $\CP^\infty$ so does their product $\C^m\times \mathcal{C}\times\Omega$. Therefore a  local Sasakian immersion of $\R^3$ with the cigar metric into $\C^m\times N\times\Omega$ would contradict \cite{LoiZedda16Cigar} where it is proven that no multiple of the cigar metric is induced by an immersion in $\CP^\infty$. 
Moreover, since $\R^3$ (with the cigar metric) is a Sasakian submanifold of $\R^{2n+1}$ (with the cigar metric), the same conclusion holds for $\R^{2n+1}$ and we have proven the following result.
\begin{prop}\label{PropSasakianCigar}
The SRS on $\R^{2n+1}$ defined by the cigar metric cannot be locally immersed in the Sasakian fiber product $\C^m\times N\times\Omega$  of the Sasakian space form $R^{2m+1}$ with a compact homogeneous Sasakian manifold $N$ and the non-compact homogeneous manifold $\Omega\times \R$ for a homogeneous bounded domain $\Omega$.
\end{prop}

Notice that in the \K setting the consequences of the presence of flat factors were already considered in \cite[Remark~1]{LoiMossa23Rigidity}.
It was proved in \cite{LoiMossa21KRSIntoComplexSpaceForms,LoiMossa22HBD,LoiMossa22BlowUp,LoiMossa23Rigidity} that non-trivial KRS cannot be locally induced by complex space forms or a certain class of homogeneous \K manifolds without flat factors.
The discussion above gives an example of KRS which cannot be induced by a homogeneous \K manifold with a flat factor as it cannot be immersed into the product $\C^m\times \mathcal{C}\times\Omega$ if $\mathcal{C}$ is of integral type.
In fact, with different and more involved techniques, we are able to generalize this to the following result.
\begin{theor}\label{ThmCigar}
	$\left(\C^n, g_{\g}\right)$ cannot be locally immersed in the product of a special generalized flag manifold  $\left(\mathcal C, g_{\mathcal C} \right)$, a  homogeneous bounded domain $\left(\Omega, g_\Omega\right)$ and a flat manifold $\left(\mathcal E, g_\me\right)$.
\end{theor}

Here we call a generalized flag manifold $\mathcal C$ \textit{special} if either its \K form is a multiple of an integral form or $\mathcal C$ is of classical type, we refer the reader to \cite{LoiMossa23Rigidity} for further details.
We believe that Theorem~\ref{ThmCigar} is an interesting result in its own right for several reasons.
Firstly, while the \K analogs of Theorems~\ref{ThmSasakiSolitonNoncompact},~\ref{ThmSasakiSolitonCompact} and \ref{ThmProducts} are established, this is new and is not implied by its Sasakian counterpart.
In fact, Theorem~\ref{ThmCigar} is stronger than Proposition~\ref{PropSasakianCigar}.
Namely, a compact homogeneous Sasakian manifold is a Boothby-Wang bundle over a generalized flag manifold $\mathcal C$. Thus $\mathcal C$ is not just special but necessarily of integral type.

\subsection*{Organization of the paper}
In the next section we review some known results about SRS and KRS. Moreover, we recall some useful constructions and present several examples. 
Section~\ref{SecKaehler} is dedicated to the proof of Theorem~\ref{ThmCigar} while in Section~\ref{SecSasaki} we prove Theorems~\ref{ThmSasakiSolitonNoncompact},~\ref{ThmSasakiSolitonCompact} and~\ref{ThmProducts}. In addition, we give an alternative proof of Proposition~\ref{PropSasakianCigar} based on Theorem~\ref{ThmCigar}.

\section{Sasakian manifolds, K\"ahler-  and Sasaki-Ricci solitons}\label{sectionbackground}

\subsection{Sasakian manifolds and immersions}

We begin with some definitions and known results, for a more exhaustive treatment we refer to the monograph by Boyer and Galicki~\cite{Boyer08Book}. 
All manifolds are assumed to be smooth, connected and oriented.

A \textit{K-contact structure} $(M,\eta,\phi,R,g)$ on a manifold $M$ consists of a contact form $\eta$ and an endomorphism $\phi$ of the tangent bundle $TM$ satisfying the following properties:
\begin{enumerate}
	\item[$\bullet$] $\phi^2=-\Id+R\otimes\eta$ where $R$ is the Reeb vector field of $\eta$,
	\item[$\bullet$] $\phi_{\vert\D}$ is an almost complex structure compatible with the symplectic form $\di\eta$ on $\D=\ker\eta$,
	\item[$\bullet$] The Reeb vector field $R$ is Killing with respect to the metric $g(\cdot,\cdot)=\di\eta(\phi\cdot,\cdot)+\eta(\cdot)\eta(\cdot)$.
\end{enumerate}   
A K-contact structure induces an almost \K structure on the metric cone $\big( S\times\R^+,t^2g+\di t^2\big)$ with the almost complex structure
\begin{enumerate}
	\item[$\bullet$] $I=\phi$ on $\D=\ker\eta$,
	\item[$\bullet$] $I(R)=t\partial_t$.
\end{enumerate} 
A \textit{Sasakian structure} is a K-contact structure $(M,\eta,\phi,R,g)$ whose metric cone is K\"ahler.
We call a manifold $M$ \textit{Sasakian} if it is equipped with a Sasakian structure. We often omit the structure $(\eta,\phi,R,g)$ itself and write simply $M$.

A Sasakian manifold is equipped with a foliation, called the \textit{Reeb foliation}, defined by the integral curves of the Reeb vector field and is called \textit{regular} (respectively \textit{quasi-regular}, \textit{irregular}) if its Reeb foliation is such.
Every regular compact Sasakian manifold $M$ is a \textit{Boothby-Wang fibration}  over a projective manifold $(X,\omega)$ with a suitable multiple $a\omega$ of $\omega$ representing an integral class (\cite{Boothby58Contact,Boyer08Book}), that is, the principal $S^1$-bundle $\pi\colon M\longrightarrow X$ with Euler class $[a\omega]$ and connection $1$-form $a\eta$ such that $\pi^*(\omega)=\di\eta$.
This is not necessarily true in the non-compact case. Nevertheless one can prove a similar statement under the regularity and completeness conditions, cf. \cite{Reinhart59Foliated}.
In general the Reeb foliation $\F$ of a Sasakian structure is transversally K\"ahler. This endows the space of leaves with a K\"ahler structure.

Given two Sasakian manifolds$(M_1,\eta_1,\phi_1,R_1,g_1)$ and  $(M_2,\eta_2,\phi_2,R_2,g_2)$, a \textit{Sasakian immersion} of $M_1$ in $M_2$ is an immersion $\psi\colon M_1\lra M_2$ such that

\begin{align*}
	\psi ^*\eta_2&=\eta_1,\hspace{20mm}\psi ^*g_2=g_1,\\
	\psi_*R_1&=R_2\hspace{5mm}\mbox{and}\hspace{5mm}\psi_*\circ \phi_1=\phi_2\circ\psi_*.
\end{align*}
Notice that this is equivalent to a holomorphic isometry $M_1\times\R^+\lra M_2\times\R^+$ of the \K cone of $M_1$ into the \K cone of $M_2$.

\subsection{K\"ahler and Sasaki-Ricci solitons}

A K\"ahler-Ricci soliton (KRS) on a complex manifold $M$ is a pair $(g, V)$ consisting of a \K metric $g$ and a holomorphic vector field $V$, called the solitonic vector field, such that
\begin{equation}\label{EqKahlerSoliton}
	\Ric_g=\lambda g+\Lie_V g 
\end{equation}
 for some $\lambda\in\R$, where $\Ric_{g}$ is the Ricci tensor of the metric $g$ and $\Lie_Vg$ denotes the
Lie derivative of $g$ with respect to $V$. Clearly, a KRS with trivial solitonic vector field is a K\"ahler-Einstein manifold. If the solitonic vector field is the gradient $\nabla f$ of a real valued function $f$, then we have a gradient KRS.

On a Sasakian manifold $(M,\eta,\phi,R,g)$ the tangent bundle splits canonically as $TM=\D\oplus T_\F$ where $\D=\ker\eta$ and $T_\F$ denotes the tangent bundle to the Reeb foliation $\F$. The transverse K\"ahler geometry is given by $(\D,\phi_{\vert_\D},\di\eta)$.
When the space of leaves of the Reeb foliation is a K\"ahler manifold $(X,J,\omega)$ we have a fibration $\pi\colon M\lra X$ such that
\begin{equation*}
	\pi^*\omega=\di\eta\hspace{5mm}\mbox{and}\hspace{5mm}\pi_*\circ\phi=J\circ\pi_*.
\end{equation*}
In virtue of this the metric decomposes as 
\begin{equation}\label{EqMetricSplit}
	g=g^T\oplus\eta\otimes\eta
\end{equation}
where $g^T(\cdot,\cdot)=\di\eta(\cdot,\phi\cdot)$.
With an abuse of notation we write $g^T$ for both the transverse metric and the metric on $X$.
It follows from \eqref{EqMetricSplit} that the Riemannian properties of $M$ can be expressed in terms of those of the transverse K\"ahler geometry and of the contact form $\eta$. For instance, the Ricci tensor of $g$ is given by
\begin{equation}\label{EqRicciSplit}
	\Ric_g=\Ric_{g^T}-2g.
\end{equation}
A Sasakian manifold $(M,\eta,\phi,R,g)$ is said to be \textit{$\eta$-Einstein} if the Ricci tensor satisfies
\begin{equation}\label{EqEtaEinstein}
	\Ric_g=\lambda g+\nu \eta\otimes\eta
\end{equation}
for some constants $\lambda,\nu \in\R$. 
It follows from \eqref{EqRicciSplit} and \eqref{EqEtaEinstein} that a Sasakian manifold is $\eta$-Einstein with constants $(\lambda,\nu)$ if, and only if, its transverse geometry is K\"ahler-Einstein with Einstein constant $\lambda+2$.

As in the \K case, Sasaki-Ricci solitons (SRS) provide a generalization of $\eta$-Einstein manifolds. we refer the reader to \cite{FutakiOnoWang09DefSRS,Placini21SasakiSolitons} for details on their definition.
A \textit{Sasaki-Ricci soliton} (SRS in short) on a Sasakian manifold $M$ is a pair $(g,V)$ consisting of the Sasakian metric $g$ and a Hamiltonian holomorphic vector field $V$ such that
\begin{equation}
	\Ric_{g^T}=\lambda g^T+\Lie_Vg^T
\end{equation}
for some $\lambda\in\R$.
If a manifold $M$ is endowed with a SRS, with an abuse of notation we will  simply say that $M$ is a SRS. 
One can easily construct examples of SRS on open Sasakian manifolds, also in the case where $\lambda\leq0$, as bundles over certain gradient KRS, see for instance \cite{Cao96ExistenceKRS,Cao97LimitsKRS}.

Notice that by definition a Sasaki-Ricci soliton $(X,g)$ on a regular Sasakian manifold is a KRS on the transverse K\"ahler geometry. In particular, if the space of leaves of the Reeb foliation is a smooth manifold $K$, then it has a canonically induced KRS. 
Viceversa, we give below some instances where one easily constructs a regular SRS from a KRS.
In order to do so, we need to recall a couple of construction of Sasakian manifolds.

Firstly observe that we can associate a regular Sasakian manifold $X\times F$ with Reeb fiber $F=\R$ or $F=S^1$ to any contractible \K manifold $(X,\omega)$. 
To see this, consider the projection $\pi : X\times F\lra X$ and let $\alpha$ be a $1$-form on $X$ such that $\di\alpha=\omega$. Define the contact form of $X\times F$ to be $\eta=\di t +\pi^*\alpha$ with Reeb vector field $\de_t$ where $t$ is the coordinate on the second factor. The endomorphism $\phi$ is defined to be the horizontal lift of the complex structure of $X$. The metric $g$ is given by $g=g^T+\eta\otimes\eta$ where $g^T$ is the metric on $X$ associated to $\omega$.
It is immediate to see that the \K metric on $X$ is a (gradient) KRS if and only if the Sasakian metric on $X\times\R$ is a (gradient) SRS.

Now consider two regular, complete Sasakian manifolds $M_1$ and $M_2$ such that they fiber on \K manifolds $X_1$ and $X_2$ respectively. 
Suppose the fibers $F$ of the fibrations are either both $\R$ or both $S^1$ so that the Reeb vector fields define a $F$-action.
Then the manifold $M$ obtained quotienting $M_1\times M_2$ by the diagonal $F$-action carries a natural Sasakian structure.
In fact, $M$ is a principal $F$-bundle over over $X_1\times X_2$ and the transverse \K structure is the product structure. 
When the manifolds are compact and $F=S^1$ one gets the so-called $(1,1)$-join $M_1*M_2$, see \cite{Boyer08Book}.
We call $M$ the \textit{Sasakian fiber product} of $M_1$ and $M_2$.

\begin{ex}\label{ExSasakianCigar}\rm
Consider on $\C$ the cigar metric
\begin{equation}\label{eqcigmetr}
	g_{\g}=\frac{d z \otimes d \bar{z}}{1+|z|^2}
\end{equation}
 due to Hamilton, see \cite{Hamilton88RicciFlow,LoiZedda16Cigar} for further details. One can see that this is a gradient KRS since $g_\g$ satisfies Equation~\eqref{EqKahlerSoliton} with $V=\nabla(\frac{1}{2}\vert z\vert^2)$.
 Since $\C$ is contractible, this induces a gradient SRS on $\R^3$ (and $\C\times S^1$ as well) which we call the \textit{Sasakian cigar metric} and denote again by $g_\g$ by abuse of notation.
 
 Notice that we one endow $\C^n$ with the product metric 
 \begin{equation}\label{eqcigmetrCn}
 	g_{\g}=\sum_{j=1}^n\frac{d z_j \otimes d \bar{z_j}}{1+|z_j|^2}
 \end{equation}
 and this is again a gradient KRS with respect to the function $f=\sum_{j=1}^n(\frac{1}{2}\vert z_j\vert^2)$.
 Therefore we can consider the gradient SRS induced on $\R^{2n+1}$. This is nothing but the Sasakian fiber product of the metrics we have just defined on $\R^3$.
\end{ex}
The next examples are obtained directly from \cite[Remark~1]{LoiMossa23Rigidity} with the constructions described above. They show that the Sasakian fiber products of certain Sasakian space forms are in fact non-trivial SRS.
\begin{ex}\label{ExSoliton1}\rm
	Consider the product metric on  $\C\times\CH^1$ where $\C$ is equipped with the flat metric and $\CH^1$ with the hyperbolic metric. This is a KRS with solitonic vector field $V=2(z\de_z+\bar{z}\de_{\bar{z}})$ where $z$ is the coordinate on $\C$.
	Since $\C\times\CH^1$ is contractible this induces a SRS on $\R^3\times\CH^1$ which is the Sasakian fiber product of the Sasakian space forms $\R^3$ and $\CH^1\times\R$.
\end{ex}

\begin{ex}\label{ExSoliton2}\rm
	Consider the product metric on  $\C\times\CP^1$ where $\C$ is equipped with the flat metric and $\CP^1$ with the Fubini-Study metric. This is a KRS with solitonic vector field $V=-2(z\de_z+\bar{z}\de_{\bar{z}})$ where $z$ is the coordinate on $\C$.
Since $\C$ is contractible there is a natural Sasakian structure on $\C\times S^1$ (the standard structure of Sasakian space form). So we can consider the Sasakian fiber product $\C\times S^3$ of the Sasakian space forms $\C\times S^1$ and $S^3$. This Sasakian structure is therefore a non-trivial SRS induced by $\C\times\CP^1$ on a manifold which is not contractible as in the previous example.
\end{ex}

\section{Local holomorphic isometries of the cigar metric}\label{SecKaehler}
Let us recall some know facts on the cigar metric $g_\g$ on $\mathbb{C}$, see \cite{LoiZedda16Cigar} for further details.
The diastasis function of the metric $g_\g$ at the origin is given by 
\begin{equation*}
	D_0^\g\left(z\right)=\int_0^{|z|^2} \frac{\log (1+t)}{t} d t.
\end{equation*}
Therefore the power series expansion of $D_0^\g\left(z\right)$ around the origin reads
\begin{equation*}
	D_0^\g\left(z\right)=\sum_{j=1}^{\infty}(-1)^{j+1} \frac{|z|^{2 j}}{j^2}.
\end{equation*}
By duplicating the variable in this last expression, we get
\begin{equation*}
	D_0^\g\left(z,\bar w\right)=\sum_{j=1}^{\infty}(-1)^{j+1} \frac{z^{ j}\bar w^j}{j^2}
\end{equation*}
so that the diastasis function $D^\g_{u_0}(z)$ centred at $u_0 \in \C$ has the following power series expansion
\begin{equation}\label{eqdcuz}
	D^\g_{u_0}(z)=\sum_{j=1}^{\infty} \frac{(-1)^{j+1}}{j^2}\left(|z|^{2 j}+|u_0|^{2 j}-z^{ j} \bar{u_0}^{ j}-u_0^{ j} \bar{z}^{ j}\right) .
\end{equation}
By duplicating the variable in this last expression, we get
\begin{equation}\label{eqdcuzw}
	D^\g_{u_0}(z,\bar w)=\sum_{j=1}^{\infty} \frac{(-1)^{j+1}}{j^2}\left(z^j\bar w^j+|u_0|^{2 j}-z^{ j} \bar{{u_0}}^{ j}-{u_0}^{ j} \bar{w}^{ j}\right).
\end{equation}

We compute here the derivative of $D^\g_{u_0}(z,\bar w)$ with respect to $\bar w$ for future reference 
\begin{equation}\label{eqderdiascig}\begin{split} 
		\frac{\de}{\de \bar w}D_{u_0}^\g\left(z,\bar w\right)&
		=\sum_{j=1}^{\infty}(-1)^{j+1} \frac{(z^j - u_0^{ j})\ \bar w^{j-1}}{j}\\
		&=\frac{\log\left(1+ z \bar w   \right)}{\bar w}-\frac{\log\left(1+u_0 \bar w\right)}{\bar w}\\
		&=\frac{1}{\bar w}\log\left(\frac{1+ z\bar w}{1+u_0 \bar w}\right)
	\end{split}
\end{equation}

Before the proof of Theorem~\ref{ThmCigar} let us recall  a useful lemma.
Let $\mathcal N_n$ be the set of real analytic functions of diastasis type $\xi : V\subset \C^n \to \R $ defined in some open domain $V\subset \C^n$, such that its real analytic extension $\xi (z, w)$ around the diagonal of  $V \times \operatorname {Conj} V$ is a holomorphic Nash algebraic function.

	\begin{lem}\label{lemtris} (\cite[Lemma 3.1]{LoiMossa23Rigidity}) and \cite[Proof of Theorem 1.1]{LoiMossa22HBD})
		Let $\mathcal C$ be a $n_1$-dimensional special generalized flag manifold and let
		$\left(\Omega, g_\Omega\right)$ be a $n_2$-dimensional homogeneous bounded domain. Then around any point $(p,q) \in \mathcal C \times \Omega$ we can find coordinates, centered at the origin, such that the diastasis function $D^{\CO}_0$ centred at the origin for the metric $g_{\mathcal C} \oplus g_\Omega $ satisfies
		\begin{equation}\label{diastomega}
			e^{D^{\CO}_0}\in \mathcal N_{n_1+n_2}^{c_{1}}\cdots \mathcal N_{n_1+n_2}^{c_{s+r}}, \qquad c_{1},\dots,c_{s+r} \in \R^+,
		\end{equation}
		where $\mathcal N_{n_1+n_2}^{c_{1}}\cdots \mathcal N_{n_1+n_2}^{c_{s+r}}
		=\left\{\xi_{1}^{c_{1}}\cdots \xi_{s+r}^{c_{s+r}}\hs\mid\hs \xi_{1},\ldots, \xi_{s+r} \in \mathcal N_{n_1+n_2} \right\}$.
	\end{lem}

	\begin{proof}[Proof of Theorem~\ref{ThmCigar}]
		Clearly it is enough to prove the statement for $n=1$ as $(\C,g_\g)$ can be embedded in $(\C^n,g_\g)$ as one of the factors.
		Let us fix a point $u_0 \in \C$. Assume by contradiction that there exist a holomorphic immersion
		\begin{align*}
	\Psi:W&\lra \left(\mc \times \Omega\right) \times \me \\
	 z&\mapsto  (\Psi_1(z),\Psi_2(z))
		\end{align*}
		of a neighbourhood   $W \subset \C$ of $u_0$ such that $g_\g=\Psi^*\left(g_{\mathcal C} \oplus g_\Omega \oplus g_\me\right) $. By the hereditary property of the diastasis function, we have
		\begin{equation}\label{eqeqdias}
			D^\g_{u_0}(z)=D^{\CO}_{\Psi_1(0)}\left(\Psi_1(z)\right)+ D^{\mathcal E}_{\Psi_2(0)}\left(\Psi_2(z)\right).
		\end{equation} 
		Passing to local coordinates centred at $\Psi_1(u_0)$ given by Lemma \ref{lemtris} and to flat coordinates centred at $\Psi_2(u_0)$, from \eqref{eqdcuz} we get
		$$
		\int_0^{|z|^2+|{u_0}|^2-z{\ov u_0} - u_0 \bar {z}} \frac{\log (1+t)}{t} d t =\sum_{k=1}^{s+r} c_k\log\left(\xi_k(\Psi_1(z)\right) + \norm{\Psi_2(z)}^2
		$$
		for suitable  $\xi_1,\ldots, \xi_{r+s} \in \mathcal N_{n_1+n_2}$. Now write
		$$
		\Psi_1(z)=\left(\psi_{1}(z), \dots, \psi_{n_1+n_2}(z)\right),
		$$
		$$
		\Psi_2(z)=\left(\psi_{n_1+n_2+1}(z), \dots, \psi_{n_1+n_2+n_3}(z)\right)
		$$
		and let $D$ be an open neighborhood of $u_0$ on which each $\psi_j$,  $j=1, \dots ,N:=n_1+n_2+n_3$, is defined. Now denote by $S$ the set  of functions $S= \left\{z,\psi_1,\ldots, \psi_{N}\right\}$.
		Consider the field  $\mathfrak R$ of rational function on $D$ and its  field extension 
		$\mathfrak F = \mathfrak R \left(S\right)$. That is, $\mathfrak F$ is the smallest subfield of the field of the  meromorphic functions on $D$, containing  rational functions and the elements of $S$. Let $l$ be the transcendence degree  of the field extension $\FF/\FR$. 
		Complete $\left\{z\right\}$ to a maximal algebraic independent subset $\left\{z,\psi_1,\ldots, \psi_l\right\} \subset S$, where $l\in \left\{0,1,\ldots,N\right\}$. With $l=0$, we mean that $\left\{z\right\}$ is already a maximal subset, i.e., each element in $S$ is holomorphic Nash algebraic.
		Then there exist minimal polynomials $P_j\left(z, X,Y\right)$ with $X=\left(X_1,\dots,X_l\right)$,  such that 
		$$
		P_j\left(z,\Phi(z), \psi_j(z)\right)\equiv 0, \ \forall j=1, \dots ,N,
		$$
		where $\Phi(z)=\left(\psi_1(z),\dots,\psi_l(z)\right)$. Moreover, by the definition of minimal polynomial,
		$$
		\frac{\de P_j\left(z,X,Y\right)}{\de Y}\left(z,\Phi (z),\psi_j(z)\right)\not\equiv 0, \ \forall j=1, \dots ,N.
		$$
		Thus, by the algebraic version of the implicit function theorem, there exist a connected open subset $U\subset D$ with $u_0\in \ov U$ and holomorphic Nash algebraic functions 
		$\hat \psi_j(z,X)$,  defined in a neighborhood $\hat U$ of $\left\{(z, \Phi(z)) \mid z \in U \right\}\subset \C \times \C^{l}$, such that
		$$
		\psi_j(z)=\hat \psi_j\left(z,\Phi(z)\right), \ \forall j=1, \dots ,N.
		$$
		for any $z\in U$.	Consider now the function 
		\begin{equation*}\begin{split} 
				F(z,X,\bar w) :&= D^\g_{u_0}(z,\bar w) - D^{\CO}_{\Psi_1(0)}\left(\hat\Psi_1(z,X),\ov{\Psi_1(w)}\right) - D^{\mathcal E}_{\Psi_2(0)}\left(\hat\Psi_2(z,X),\ov{\Psi_2(w)}\right)\\
		\end{split}\end{equation*}
		where 
		$$
		\hat\Psi_1(z,X):=\left(\hat\psi_{1}(z,X), \dots,\hat \psi_{n_1+n_2}(z,X)\right)
		$$
		and
		$$
		\hat\Psi_2(z,X):=\left(\hat\psi_{n_1+n_2+1}(z,X), \dots,\hat \psi_{N}(z,X)\right).
		$$ 
		It is not restrictive to assume that $F$ is defined on $\hat U \times D$. We claim that 
		$$
		\left(\frac{\de}{\de \bar w}\left(\bar w\hs\frac{\de}{\de \bar w} F\right)\right)(z,X,\bar w)\equiv 0
		$$
		for all $w \in U$. Assume, by contradiction, that there exists $w_0\in U$ such that 
		$
		\left(\frac{\de}{\de \bar w}\left(\bar w \frac{\de}{\de \bar w} F\right)\right)(z,X,\bar w_0)\neq 0.
		$
		Making use of \eqref{eqderdiascig} we get
		\begin{equation}\label{eqdewF}\begin{split} 
				\frac{\de}{\de \bar w} F(z,X,\bar w)&=\frac{1}{\bar w}\log\left(\frac{1+ z\bar w}{1+u_0 \bar w}\right)\\
				&-\sum_{k=1}^{r+s} c_k\frac{
					\frac{\de }{\de \bar w}
					\left(\xi_k\left(\hat\Psi_1(z,X),\ov {\Psi_1(w)}\right)\right)
				}{\xi_k\left(\hat\Psi_1(z,X),\ov {\Psi_1(w)}\right)}
				- \left\langle \hat\Psi_2(z,X), \ov{\Psi_2'(w)} \right\rangle.
		\end{split}\end{equation}
	
		Thus, we see that $\left(\frac{\de}{\de \bar w}\left(\bar w\frac{\de}{\de \bar w} F\right)\right)(z,X,\bar w_0)$ is Nash algebraic in $(z,X)$.
		Hence there exists a holomorphic polynomial $P(z,X,t)=A_d(z,X)t^d+\dots+A_0(z,X)$ with $A_0(z,X)\not\equiv 0$ such that  
		$$
		P\left(z,X,\left(\frac{\de}{\de \bar w}\left(\bar w\frac{\de}{\de \bar w} F\right)\right)(z,X,\bar w_0)\right)=0.
		$$ 
		Since  by construction $F(z,\Phi (z),\bar w)\equiv 0$, we see that  
		$
		\left(\frac{\de}{\de \bar w}\left(\bar w\frac{\de}{\de \bar w} F\right)\right)(z,\Phi (z),\bar w)\equiv 0.
		$
		Thus  $A_0(z,\Phi (z))\equiv 0$,
		which contradicts  the fact that $\psi_1(z),\dots,\psi_l(z)$ are algebraically independent over $\FR$. Hence $\left(\frac{\de}{\de \bar w}\left(\bar w \frac{\de}{\de \bar w} F\right)\right)(z,X,\bar w_0)\equiv 0$ and the claim is proved. 
		
		We conclude that $ \frac{\de}{\de \bar w} F(z,X,\bar w)$ satisfies the following differential equation
		$$
		\frac{\de}{\de \bar w} F(z,X,\bar w) + \bar w \frac{\de^2}{\de \bar w^2} F(z,X,\bar w)\equiv 0.
		$$
This implies
		\begin{equation}\label{eqsolcdiff}
			\frac{\de}{\de \bar w} F (z,X,\bar w)= \frac{c(z,X)}{\bar w} 
		\end{equation}
	where, using  \eqref{eqdewF} and the fact that $c(z,X)$ is independent of $\bar w$, we see that
		\begin{equation}\begin{split} \label{eqexpc}
				c(z,X)  =& \bar u_0 \frac{\de F}{\de \bar w} (z,X,u_0)=\log\left(\frac{1+ z\bar u_0}{1+\abs {u_0} ^2}\right)\\
				- & \bar u_0\sum_{k=1}^{r+s}c_k \frac{
					\frac{\de }{\de \bar w}
					\left(\xi_k\left(\hat\Psi_1(z,X),\ov {\Psi_1(u_0)}\right)\right)
				}{\xi_k\left(\hat\Psi_1(z,X),\ov {\Psi_1(u_0)}\right)}- \bar u_0 \left\langle \hat\Psi_2(z,X), \ov{\Psi_2'(u_0)} \right\rangle.
		\end{split}\end{equation}

	 Combining  \eqref{eqdewF}, \eqref{eqsolcdiff} and \eqref{eqexpc}  we get
		$$
		\log\left(\frac{1+ z\bar w}{1+u_0 \bar w}\right)
		-\bar w\sum_{k=1}^{r+s}c_k \frac{
			\frac{\de }{\de \bar w}
			\left(\xi_k\left(\hat\Psi_1(z,X),\ov {\Psi_1(w)}\right)\right)
		}{\xi_k\left(\hat\Psi_1(z,X),\ov {\Psi_1(w)}\right)}-  \bar w \left\langle \hat\Psi_2(z,X), \ov{\Psi_2'(w)} \right\rangle
		$$
		$$
		= \log\left(\frac{1+ z\bar u_0}{1+\abs {u_0} ^2}\right) 
		-  \bar u_0 \sum_{k=1}^{r+s}c_k \frac{
			\frac{\de }{\de \bar w}
			\left(\xi_k\left(\hat\Psi_1(z,X),\ov {\Psi_1(u_0)}\right)\right)
		}{\xi_k\left(\hat\Psi_1(z,X),\ov {\Psi_1(u_0)}\right)}- \bar u_0 \left\langle \hat\Psi_2(z,X), \ov{\Psi_2'(u_0)} \right\rangle.
		$$
		Hence
		$$
		\log\left(\frac{1+ z\bar w}{1+ z\bar u_0}\cdot\frac{1+\abs {u_0} ^2}{1+u_0 \bar w}\right)=
		$$
		$$
		\bar w\sum_{k=1}^{r+s} c_k\frac{
			\frac{\de }{\de \bar w}
			\left(\xi_k\left(\hat\Psi_1(z,X),\ov {\Psi_1(w)}\right)\right)
		}{\xi_k\left(\hat\Psi_1(z,X),\ov {\Psi_1(w)}\right)}
		+  \bar w \left\langle \hat\Psi_2(z,X), \ov{\Psi_2'(w)} \right\rangle
		$$
		$$
		- \bar u_0 \sum_{k=1}^{r+s} c_k\frac{
			\frac{\de }{\de \bar w}
			\left(\xi_k\left(\hat\Psi_1(z,X),\ov {\Psi_1(u_0)}\right)\right)
		}{\xi_k\left(\hat\Psi_1(z,X),\ov {\Psi_1(u_0)}\right)}
		- \bar u_0 \left\langle \hat\Psi_2(z,X), \ov{\Psi_2'(u_0)} \right\rangle,
		$$
		where the right hand side of the equality is Nash algebraic with respect to $(z,X)$. Thus,
		applying \cite[Lemma 2.2]{HuangYuan15SubmanifoldsHSS}, we deduce that $\log\left(\frac{1+ z\bar w}{1+ z\bar u_0}\cdot\frac{1+\abs {u_0} ^2}{1+u_0 \bar w}\right)$
		is constant in $z$ for any fixed $w \in D$, which gives us the desired contradiction. The proof is complete.
	\end{proof}

\section{Local Sasakian immersion of Sasaki-Ricci solitons}\label{SecSasaki}
In this section we collect a discussion on Sasakian immersions of SRS into homogeneous Sasakian manifolds and prove the theorems stated in the introduction.
Firstly we use Sasakian rigidity \cite[Theorem~1.2]{LoiPlaciniZedda23SasakiHomogeneous} to drop the regularity assumption in Theorem~1 and Theorem~2 of \cite{Placini21SasakiSolitons} and obtain the following.

\begin{proof}[Proof of Theorem~\ref{ThmSasakiSolitonCompact} and Theorem~\ref{ThmSasakiSolitonNoncompact}]
	Let $M$ be a complete Sasakian manifold endowed with a Sasaki-Ricci soliton. 
		Suppose there exist a neighbourhood $U_p$ of a point $p\in M$ and an immersion $\psi\colon U_p\lra S(N,c)$ into a Sasakian space form $S(N,c)$ with $c\in\R$.
	
	Firstly, consider the universal cover $\widetilde{M}$ with the induced Sasakian-Ricci soliton. 
	Notice that the metric on $\widetilde{M}$ is real analytic because it is a SRS, see \cite[Corollary~1.3]{Kotschwar13SolitonRealAnalytic}. 
	Therefore, by a direct  application of \cite[Theorem~8]{Calabi53Isometric} in the Sasakian setting, for every point  $q\in \widetilde{M}$ there exists a neighbourhood $U_q$ that admits an embedding  $U_q\lra S(N,c)$ in the Sasakian space form $S(N,c)$.
	This, together with $\pi_1(\widetilde{M})=0$, shows that $\widetilde{M}$ satisfies the hypotheses of \cite[Theorem~1.2]{LoiPlaciniZedda23SasakiHomogeneous} so there exists a global Sasakian immersion of $\widetilde{M}$ into a Sasakian space form. 
	It follows from \cite[Theorem~1]{Placini21SasakiSolitons} that $\widetilde{M}$ is a Sasakian space form $S(n,c)$ if $c\leq -3$. 
	Correspondingly, if $c=1$ (i.e. if $S(N,c)$ is the standard Sasakian sphere $S^{2N+1}$),  then  \cite[Theorem~2]{Placini21SasakiSolitons} implies that $\widetilde{M}$ is $\eta$-Einstein with $\eta$-Einstein constants $(\lambda,\nu)$ given by $\lambda=4\mu -2$ for some $\mu\in\Q$. 
	The theses then follow from the fact that $M$ is the quotient of $\widetilde{M}$ by a discrete group of Sasakian automorphisms.
\end{proof}

\begin{proof}[Proof of Corollary~\ref{CorHomogeneous}]
	It was proven in \cite[Theorem~1.4]{LoiPlaciniZedda23SasakiHomogeneous} that, up to trasverse homotheties, every compact homogeneous Sasakian manifold admits a Sasakian embedding into a Sasakian space form.
	Thus, after possibly deforming the SRS by a transverse homothety, we get a immersion of $\varphi:U\lra S^{2N+1}$ in standard Sasakian sphere.
	We deduce that the transversally rescaled SRS is trivial, that is $\eta$-Einstein, and consequently so is the starting Sasakian metric.
\end{proof}

One can easily construct Sasakian manifolds as Sasakian fiber products of Sasakian space forms (or, more generally, regular Sasakian manifolds) by identifying the leaves of the Reeb foliation, see Section~\ref{sectionbackground} above.
In analogy with the \K case, one can ask whether SRS satisfy similar constraints when they can be locally  immersed into fiber products of Sasakian space forms. Example~\ref{ExSoliton1} and Example~\ref{ExSoliton2} above shows that this is generally false.
Nevertheless, in light of the recent theorem of the first named author and Loi \cite[Theorem~1.1]{LoiMossa23Rigidity}, we derive the following result.

\begin{proof}[Proof of Theorem~\ref{ThmProducts}]
	Since the immersion is local and SRS lift to the covering spaces, we can replace $M$ by its universal cover, that is, we can assume that $\pi_1(M)=0$.
	Now we show, similarly to the proof of Theorems~\ref{ThmSasakiSolitonNoncompact} and~\ref{ThmSasakiSolitonCompact}, that the Sasakian structure on $M$ is regular.
	Notice that a suitable transverse homothety of the Sasakian manifold $M'$ can be immersed into the infinite dimensional Sasakian space form $S^\infty$ by \cite[Theorem~1.5]{LoiPlaciniZedda23SasakiHomogeneous}.
	Therefore, the same $\D$-homothety of $M$ admits a local embedding $\psi:U_p\lra S^\infty$ and is again a SRS.
	Moreover, the metric on $M$ is real analytic because it is a SRS, see \cite[Corollary~1.3]{Kotschwar13SolitonRealAnalytic}. 
	Thus, by a direct  application of \cite[Theorem~8]{Calabi53Isometric} in the Sasakian setting, for every point  $q\in M$ there exists a neighbourhood $U_q$ that admits an embedding  $U_q\lra S^\infty$.
	This, together with $\pi_1(M)=0$, shows that $M$ satisfies the hypotheses of \cite[Theorem~1.2]{LoiPlaciniZedda23SasakiHomogeneous} so there exists a global Sasakian immersion of $\widetilde{M}$ into a Sasakian space form. 
	This in particular implies that $M$ is regular as claimed.
	Thus, by \cite{Reinhart59Foliated} $M$ admits a Boothby-Wang fibration $\pi:M\lra X$ over a KRS $X$ being a complete and regular. 
	Moreover, the immersion $\psi$ covers a holomorphic isometry $\varphi$ of the neighbourhood $V_x=\pi(U_p)$ of the point $x=\pi(p)$ into $\Omega\times F$ where $F$ is a generalized flag manifold. Visually,
	$$
	\begin{tikzcd}
		U_p \arrow[r, "\psi"] \arrow[d, "\pi"'] & \Omega\times N  \arrow[d, "\pi'"] \\
		V_x \arrow[r, "\phi"'] &   \Omega\times F
	\end{tikzcd}
	$$
	Since the manifold $F$ is the base of a compact homogeneous Sasakian manifold, it is necessarily of integral type (i.e. a suitable multiple of its \K class is integral.
	Now \cite[Theorem~1.1]{LoiMossa23Rigidity} ensures that the KRS on $X$ is trivial.
	Thus $M$ (and all its $\D$-homotheties) is an $\eta$-Einstein manifold.
\end{proof}
We conclude this section giving an alternative proof of Proposition~\ref{PropSasakianCigar} using Theorem~\ref{ThmCigar}.
\begin{proof}[Proof of Proposition~\ref{PropSasakianCigar}]
	We want to prove this by contradiction. Suppose there exists a local Sasakian immersion $\varphi: U_p\lra \C^{m}\times N\times\Omega$  of a neighbourhood $U_p$ of a point $p=(t,z_1,\ldots,z_n)$.
	Since the statement is local and a compact homogeneous Sasakian manifold embeds into a sphere up to $\D$-homothety by \cite[Theorem~1.4]{LoiPlaciniZedda23SasakiHomogeneous}, we can replace $N$ by a suitable transverse homothety $D_a$ of the standard sphere $S^{2k+1}$.
	That is, the immersion $\varphi$ is of the form $\varphi: U_p\lra \C^{m}\times S_a^{2k+1}\times\Omega$ where $S_a^{2k+1}$ denotes the $\D_a$-homothety of standard sphere.
	Now this immersion covers a holomorphic isometry of a neighbourhood $V_z$ of $z=(z_1,\ldots,z_n)\in\C^n$ endowed with the product of cigar metrics into $\C^{m}\times \CP_a^{k}\times\Omega$ where $\CP_a^{k}$ is the complex projective space endowed with the \K metric $a g_{FS}$. 
	Clearly, $\CP_a^{m}$ is a flag manifold of integral type and therefore the existence of an holomorphic isometry $V_z\lra \C^{m}\times \CP_a^{k}\times\Omega$ violates Theorem~\ref{ThmCigar}.
\end{proof}

\bibliographystyle{amsplain}

\bibliography{/Users/giovanniplacini/Library/CloudStorage/Dropbox/Projects/RTT_PROGETTI/Bibliography/biblio.bib}

\providecommand{\bysame}{\leavevmode\hbox to3em{\hrulefill}\thinspace}
\providecommand{\MR}{\relax\ifhmode\unskip\space\fi MR }
\providecommand{\MRhref}[2]{%
  \href{http://www.ams.org/mathscinet-getitem?mr=#1}{#2}
}
\providecommand{\href}[2]{#2}
\begin{thebibliography}{10}

\bibitem{Boothby58Contact}
W.~M. Boothby and H.~C. Wang, \emph{On contact manifolds}, Ann. of Math. (2)
  \textbf{68} (1958), 721--734. \MR{112160}

\bibitem{Boyer08Book}
Charles~P. Boyer and Krzysztof Galicki, \emph{Sasakian geometry}, Oxford
  Mathematical Monographs, Oxford University Press, Oxford, 2008. \MR{2382957}

\bibitem{Calabi53Isometric}
Eugenio Calabi, \emph{Isometric imbedding of complex manifolds}, Ann. of Math.
  (2) \textbf{58} (1953), 1--23. \MR{57000}

\bibitem{Cao96ExistenceKRS}
Huai-Dong Cao, \emph{Existence of gradient {K}\"{a}hler-{R}icci solitons},
  Elliptic and parabolic methods in geometry ({M}inneapolis, {MN}, 1994), A K
  Peters, Wellesley, MA, 1996, pp.~1--16. \MR{1417944}

\bibitem{Cao97LimitsKRS}
\bysame, \emph{Limits of solutions to the {K}\"{a}hler-{R}icci flow}, J.
  Differential Geom. \textbf{45} (1997), no.~2, 257--272. \MR{1449972}

\bibitem{Collins16SRFlowConvergenceCRELLE}
Tristan~C. Collins, \emph{Stability and convergence of the {S}asaki-{R}icci
  flow}, J. Reine Angew. Math. \textbf{716} (2016), 1--27. \MR{3518371}

\bibitem{Collins15SRFlowConvergence}
Tristan~C. Collins and Adam Jacob, \emph{On the convergence of the
  {S}asaki-{R}icci flow}, Analysis, complex geometry, and mathematical physics:
  in honor of {D}uong {H}. {P}hong, Contemp. Math., vol. 644, Amer. Math. Soc.,
  Providence, RI, 2015, pp.~11--21. \MR{3372457}

\bibitem{FutakiOnoWang09DefSRS}
Akito Futaki, Hajime Ono, and Guofang Wang, \emph{Transverse {K}\"{a}hler
  geometry of {S}asaki manifolds and toric {S}asaki-{E}instein manifolds}, J.
  Differential Geom. \textbf{83} (2009), no.~3, 585--635. \MR{2581358}

\bibitem{Hamilton88RicciFlow}
Richard~S. Hamilton, \emph{The {R}icci flow on surfaces}, Mathematics and
  general relativity ({S}anta {C}ruz, {CA}, 1986), Contemp. Math., vol.~71,
  Amer. Math. Soc., Providence, RI, 1988, pp.~237--262. \MR{954419}

\bibitem{HuangYuan15SubmanifoldsHSS}
Xiaojun Huang and Yuan Yuan, \emph{Submanifolds of {H}ermitian symmetric
  spaces}, Analysis and geometry, Springer Proc. Math. Stat., vol. 127,
  Springer, Cham, 2015, pp.~197--206. \MR{3445521}

\bibitem{Kotschwar13SolitonRealAnalytic}
Brett~L. Kotschwar, \emph{A local version of {B}ando's theorem on the
  real-analyticity of solutions to the {R}icci flow}, Bull. Lond. Math. Soc.
  \textbf{45} (2013), no.~1, 153--158. \MR{3033963}

\bibitem{LoiPlaciniZedda23SasakiHomogeneous}
A.~Loi, G.~Placini, and M.~Zedda, \emph{Immersions into {S}asakian space
  forms}, 9 May 2023, preprint arXiv:2305.05509 [math.DG].

\bibitem{LoiMossa21KRSIntoComplexSpaceForms}
Andrea Loi and Roberto Mossa, \emph{K\"{a}hler immersions of
  {K}\"{a}hler-{R}icci solitons into definite or indefinite complex space
  forms}, Proc. Amer. Math. Soc. \textbf{149} (2021), no.~11, 4931--4941.
  \MR{4310116}

\bibitem{LoiMossa22HBD}
\bysame, \emph{Holomorphic isometries into homogeneous bounded domains}, Proc.
  Amer. Math. Soc. \textbf{151} (2023), no.~9, 3975--3984. \MR{4607641}

\bibitem{LoiMossa22BlowUp}
\bysame, \emph{On {H}olomorphic {I}sometries into {B}low-{U}ps of {${\mathbb
  {C}}^n$}}, Mediterr. J. Math. \textbf{20} (2023), no.~4, 230. \MR{4597721}

\bibitem{LoiMossa23Rigidity}
\bysame, \emph{Rigidity properties of holomorphic isometries into homogeneous
  {K}\"ahler manifolds}, 9 May 2023, preprint arXiv:2305.05244 [math.DG].

\bibitem{LoiZedda16Cigar}
Andrea Loi and Michela Zedda, \emph{On {C}alabi's diastasis function of the
  {C}igar metric}, J. Geom. Phys. \textbf{110} (2016), 269--276. \MR{3566114}

\bibitem{Petrecca16SRSDeformations}
David Petrecca, \emph{On sasaki–ricci solitons and their deformations},
  Advances in Geometry \textbf{16} (2016), no.~1, 57--70.

\bibitem{Placini21SasakiSolitons}
Giovanni Placini, \emph{Sasakian immersions of {S}asaki-{R}icci solitons into
  {S}asakian space forms}, J. Geom. Phys. \textbf{166} (2021), Paper No.
  104265, 7. \MR{4249127}

\bibitem{Reinhart59Foliated}
Bruce~L. Reinhart, \emph{Foliated manifolds with bundle-like metrics}, Ann. of
  Math. (2) \textbf{69} (1959), 119--132. \MR{107279}

\bibitem{Smoczyk10SRFlowDefinition}
Knut Smoczyk, Guofang Wang, and Yongbing Zhang, \emph{The {S}asaki-{R}icci
  flow}, Internat. J. Math. \textbf{21} (2010), no.~7, 951--969. \MR{2671532}

\bibitem{Tadano15GapSRS}
Homare Tadano, \emph{Gap theorems for compact gradient {S}asaki-{R}icci
  solitons}, Internat. J. Math. \textbf{26} (2015), no.~4, 1540009, 17.
  \MR{3338073}

\end{thebibliography}

\end{document}